\documentclass[11pt,a4paper]{article}
\usepackage{amsfonts,amsbsy,amssymb,amsmath,graphicx}
\usepackage{epsfig}
\usepackage{dsfont}
\usepackage{mathrsfs}
\usepackage{geometry}
\usepackage[T1]{fontenc}
\usepackage{subeqnarray}
\usepackage{enumerate}
\usepackage{color}
\usepackage{rotating}
\usepackage{epsfig}
\usepackage{ifpdf}
\usepackage{cite}
\usepackage[colorlinks=true,linkcolor=blue,citecolor=blue]{hyperref}

\newenvironment{proof} {\noindent {\em \textbf{Proof}} } { \hfill \fbox{~} \\ }

\def\1{1\kern-.20em {\rm l}}
\def\abs  #1\par{\vskip2truemm   {\noindent} {\parindent=12truemm\narrower\baselineskip=3truemm \pf #1\par}}
\newtheorem{theorem}{Theorem}[section]

\newtheorem{definition}[theorem]{Definition}

\numberwithin{equation}{section}

\font\bfit = cmbxti9
\newcommand{\N}{\mathbb{N}}

\newcommand{\C}{\mathbb{C}}

\def\1{1\kern-.20em {\rm l}}
\def\qed{\ifmmode\mbox{\hfill\sqb}\else{\ifhmode\unskip\fi%
\nobreak\hfil
\penalty50\hskip1em\null\nobreak\hfil\sqb
\parfillskip=0pt\finalhyphendemerits=0\endgraf}\fi}
\def\cqfd{\ifmmode\sqw\else{\ifhmode\unskip\fi\nobreak\hfil
\penalty50\hskip1em\null\nobreak\hfil\sqw
\parfillskip=0pt\finalhyphendemerits=0\endgraf}\fi}
\title{\bf The symmetric Dunkl-classical  orthogonal polynomials revisited }% The symmetric Dunkl OPS revisited}
\author{\small Khalfa DOUAK\footnote{Email: khalfa.douak@gmail.com}\\
\small Laboratoire Jacques-Louis Lions, Sorbonne Universit\'{e},\\ 75252 Paris Cedex 05, France\\
\bfit In memory of Pascal MARONI}
\date{\today}
\begin{document}
\maketitle
\begin{abstract}
\noindent We investigate  the symmetric Dunkl-classical orthogonal polynomials by using a new approach applied in connection with the Dunkl operator.
 The main aim of this technique is to determine the recurrence coefficients first and foremost.
We establish the existence and uniqueness of symmetric Dunkl-classical orthogonal polynomials, and also confirm that only two families of orthogonal polynomials, that is, the generalized Hermite polynomials and the generalized Gegenbauer polynomials, belong to this class.
Two apparently new characterizations in a more general setting are given. This paper complements earlier work of Ben Cheikh and Gaied.
\end{abstract}
%%%%%%%%%%%%%%%%%%%%%%%%%%%%%%%%%%%%%%%%
\section{Introduction and motivation}
The original Dunkl operators as called in the literature are a family of first-order  differential-difference operators first introduced by Charles Dunkl
in \cite{Dunk} whose the main idea was to use finite reflection groups to provide structure for multi-variable analysis. Here we are concerned with the
 one-dimensional Dunkl operator, denoted  by $T_\mu$,  defined as follows
 \begin{align}
 \big(T_{\mu}f\big)(x)  = f'(x)+\mu \frac{f(x)-f(-x)}{x}, \label{eq1.1}
 \end{align}
where $f$ is a polynomial or an entire function and  $\mu$ is taken here to be an arbitrary real parameter. \\
Acting on polynomials, one sees at once that the Dunkl operator $T_\mu$ reduces the degree of every polynomial by exactly one.
 In particular, if we take  $f(x)=x^n$, we get $T_\mu x^n=\mu_n x^{n-1}$ where   $\mu_n=n+\mu\left(1-(-1)^n\right),\ n\geqslant0$,
  with  $\mu\ne-n-\frac{1}{2},\ n\geqslant0$. \\
We are interested here in the problem of finding  the sequences of symmetric orthogonal polynomials  $\{P_n\}_{n\geqslant0}$ when the sequence
of polynomials $\{T_\mu P_n\}_{n\geqslant0}$ is also orthogonal.\\
In an earlier work \cite{BeGa}, Ben Cheikh and Gaied considered this problem and referred to the resulting polynomials as the symmetric
{\it Dunkl-classical} orthogonal polynomials.
In doing so, the authors have focused their study on the identification of  solutions of a differential-difference equation of type
 \begin{align}
\sigma(x)T_{\mu}^2P_n(x) +\tau(x)T_{\mu}P_n(x)+\lambda_nP_n(x) = 0,\label{eq1.2}
 \end{align}
 where  $\sigma(x)$ and $\tau(x)$ are two polynomials independent to $n$, $\sigma(x)$ is an even polynomial
 of degree less than or equal to two, $\sigma(x)$
 is an odd polynomial of  first degree and $\lambda_n$ is independent of $x$. They deduced that the generalized Hermite polynomials and
  the generalized Gegenbauer  polynomials are the only eigenfunctions of \eqref{eq1.2}. \par
 \medskip
\noindent In this paper we investigate this kind of polynomials by adopting  a different approach.
 This technique first used in my Ph.D. thesis \cite{Doua} to construct the classical orthogonal polynomials, part of which was published in \cite{DoMa1}.
 It was also applied later in many other works, for instance,  to construct the families of $D_\omega$-classical  orthogonal polynomials \cite{Doua1}.
 To continue in this line, it is then reasonable to ask whether this can be extended to the Dunkl context, where the ordinary derivative $D$
 (or $D_\omega$)  is replaced by the Dunkl operator $T_\mu$. This allows us to construct the symmetric Dunkl-classical orthogonal polynomials.
 Unfortunately, this method is rather complicated, if not impossible, in the general case. Therefore our interest is in describing such orthogonal polynomials
only  for the symmetric case.\\
Roughly speaking, the main idea behind this approach is to focus first on the explicit determination of recurrence coefficients.
This in fact allows us to establish a link between the Dunkl-classical character and the  form of the recurrence coefficients which represents
a new result characterizing    the symmetric   Dunkl-classical orthogonal polynomials (Theorem 2.1).
As an important consequence of this theorem we also bring another result characterizing the resulting polynomials
via a differential-difference equation in a more general setting than that stated in \eqref{eq1.2} (Theorem 2.2).\\
  As an application of Theorem 2.1, we can assert that the generalized Hermite polynomials and the generalized Gegenbauer polynomials are the only
  symmetric Dunkl-classical orthogonal polynomials. Subsequently, we sharpen the results stated in  Theorem 2.2 to characterize each of the above two
  families of polynomials by means of a differential-difference equation of type \eqref{eq1.2}.
 Before we can state the main theorems   in this paper we wish to introduce some standard notations.\\
 \noindent
Let $\mathscr P$ be the vector space of polynomials of one variable with complex coefficients and $\mathscr P'$ its algebraic dual.
 By $\{P_n\}_{n\geqslant0}$ we denote a polynomials sequence (PS in short), with $\deg P_n=n$, and $\{u_n\}_{n\geqslant0}$
 its associated dual sequence defined by
 $\bigl< u_n , P_m \bigr> = \delta_{n m} ;\ n , m\geqslant0,$ where $\delta_{n m}$ is the Kronecker's delta symbol.
 Throughout this article, we will always  consider  the sequence of {\it monic} polynomials, i.e. the leading coefficient of each polynomial $P_n$
 is one ($P_n(x)=x^n+\cdots$).
 We also recall the two linear operators:\\
 The ordinary differential operator $ D=d/dx$ : $\left(Df\right)(x)=f'(x)$.\\
 The Hahn operator $\, H_{q,\omega}$ defined  for two fixed complex numbers $q$ and $\omega$ by
 \begin{equation}
\left(H_{q,\omega}f\right)(x):=\frac{f(qx+\omega)-f(x)}{(q-1)x+\omega},\quad \forall f \in \mathscr P. \label{eq1.3}
\end{equation}
When $q=1$ with $\omega\ne0$ we get the discrete operator $\left(D_{\omega}f\right)(x):=\left(H_{1,\omega}f\right)(x)$, that is,
\begin{equation}
\left(D_{\omega}f\right)(x)=\frac{f(x+\omega)-f(x)}{\omega}.\label{eq1.4}
\end{equation}
For $q\ne1$ and $\omega=0$, the $q$-difference operator, which we write
$\left({H}_{q}f\right)(x):=\left(H_{q,0}f\right)(x)$, is   given by
\begin{equation}
\left(H_{q}f\right)(x):=\frac{f(qx)-f(x)}{(q-1)x}.\label{eq1.5}
\end{equation}
 With this notation, we can write \eqref{eq1.1} in the form
 \begin{align}
 \big(T_{\mu}f\big)(x)  = \left(Df\right)(x)+2\mu \left(H_{-1}f\right)(x).\label{eq1.6}
 \end{align}
To begin with, we consider a sequence of monic orthogonal polynomials (MOPS) with its  associated linear functional denoted  $u\in \mathscr P'$.
The sequence of complex numbers $(u)_n,\ n=0, 1, 2, \ldots,$ denotes the moments of $u$ with respect to the sequence $\{x^n\}_{n \geqslant0}$, namely,
the moment of order $n$ for the functional $u$ is denoted by $(u)_n := \big< u , x^n \big>, \ n=0, 1, 2, \ldots$.
Thus, the linear functional $u$ is completely determined by its moments. To terminate we define by $\Delta_n=\det\big((u)_{i+j}\big)_{i,j=0}^n,\ n\geqslant0$,
 the ordinary Hankel determinants constructed from the moments of $u$.
  \begin{definition}
A {\rm PS} $\{P_n\}_{n\geqslant0}$ is said to be orthogonal with respect to (w.r.t.) the linear functional $u$, if it satisfies the orthogonality conditions
\begin{subequations}
\begin{align}
\big< u \, ,\, P_n P_m \big> &= 0,\ n\ne m,\label{eq1.7a}\\
\big< u \, ,\, P_n^2 \big> &\not= 0,\ n\geqslant0.\label{eq1.7b}
\end{align}
 \end{subequations}
\end{definition}
In this case \eqref{eq1.7b} are said to be the regularity conditions. As an immediate consequence of the regularity of $u$, we have  $(u)_0 \ne0$
 and $u=\lambda u_0$ with $\lambda\ne0$. For this reason, we always consider the orthogonality w.r.t. $u_0$ which is the {\it canonical} form associated
 to the PS $\{P_n\}_{n\geqslant0}$.
\begin{theorem}{\rm {\rm \cite{Chih2}}}
 A necessary and sufficient condition for the existence of an orthogonal polynomials sequence {\rm (OPS)} $\{P_n\}_{n\geqslant0}$ with respect to  $u_0$
 is $\Delta_n=\det\big((u_0)_{i+j}\big)_{i,j=0}^n\ne0$, for $n=0, 1, 2, \ldots$. The functional $u_0$ is then said to be regular.
\end{theorem}
Equivalently, the polynomials $P_n,\ n=0, 1, 2, \ldots$, are orthogonal if and only if they satisfy the well-known second-order recurrence relation
    \begin{subequations}
\begin{align}
&P_{n+2}(x)=(x-\beta_{n+1})P_{n+1}(x)-\gamma_{n+1}P_n(x),\ n\geqslant0,\label{eq1.8a}\\
& P_0(x) = 1,\ P_1(x) = x-\beta_0,\label{eq1.8b}
\end{align}
    \end{subequations}
with the regularity conditions $\gamma_{n+1}\ne0,\ n\geqslant0$.\\
The linear functional $u$ is said to be regular if there exists a sequence of polynomials $\{P_n\}_{n\geqslant0}$, such that\\
{\bf 1.}  $\deg P_n=n$, $\forall n\in \N$,\\
{\bf 2.} $\big< u \, ,\, P_n P_m \big>= k_n\delta_{nm},\  k_n\ne0,\ \forall n, m\in \N$.
If  $k_n>0$, for all $n\in\N$, then $u$ is positive definite.\par
\noindent Note that  we work here with  linear functionals which we assume  to be  regular, but not necessarily positive definite.\\
A regular linear functional $u$ is said to be semi-classical if there exist two polynomials $\phi$ and $\psi$, with $\deg\phi\geqslant0$ and
$\deg\psi\geqslant1$, satisfying \cite{Maro}
\begin{align}
D\left(\phi(x)u\right)+\psi(x)u=0.\label{eq1.9}
\end{align}
If moreover $\phi(x)=a_tx^t+...$, and $\psi(x)=b_px^p+...$,
then one of the following conditions holds:\\
{$\bf a$.}  $t-1\ne p$,\\
{$\bf b$.} if $t-1= p$, then  $na_t\ne b_p, \forall n\in \N$, and the class of $u$ is defined as
\begin{align}
s=\max\left(\deg\phi-2,\deg\psi-1\right).\label{eq1.10}
\end{align}
Under these conditions, the pair $(\phi,\psi)$ is said admissible.\par\medskip
\noindent  Define now the sequence $\{Q_n\}_{n\geqslant0}$, where  $Q_n(x)= \mu_{n+1}^{-1}\big(T_\mu P_{n+1}\big)(x),\, n\geqslant0$.
This normalization plays an important role  in what follows and its choice seems to be the best adapted to our approach.
 The classical orthogonal polynomial sequences of Hermite, Laguerre, Jacobi and Bessel collectively satisfy
the so-called Hahn property \cite{Hahn1}: the sequence of its derivatives is again an orthogonal polynomial sequence.
 Likewise, we can introduce the adjective {\it Dunkl-classical} or simply {\it ``$T_\mu$-classical”} of an OPS as follows.
 \begin{definition}
The {OPS} $\{P_n\}_{n\geqslant0}$ is called ``$T_\mu$-classical” if the sequence  $\{Q_n\}_{n\geqslant 0}$ is also a {OPS}.
\end{definition}
The  $\{P_n\}_{n\geqslant0}$ is said to be {symmetric} if each polynomial fulfills $ P_n(- x)=(-1)^n P_n(x),\, n\geqslant0$.
Further, when $\{P_n\}_{n\geqslant0}$ is an OPS, all of the moments of the odd order of the form $u_0$ are zero.
\begin{theorem}{\rm \cite{Chih2}}
 For every {\rm MOPS} $\{P_n\}_{n\geqslant0}$,  the following statements are equivalent:\\
{\rm (1)} The form $u_0$ is symmetric.\\
{\rm (2)} The sequence $\{P_n\}_{n\geqslant0}$ is symmetric.\\
{\rm (3)} The sequence $\{P_n\}_{n\geqslant0}$ satisfies the recurrence relation
   \begin{subequations}
\begin{align}
&P_{n+2}(x)=xP_{n+1}(x)-\gamma_{n+1}P_n(x),\ n\geqslant0,\label{eq1.11a}\\
& P_0(x) = 1,\ P_1(x) = x.\label{eq1.11b}
\end{align}
    \end{subequations}
 \end{theorem}
 The following identity is straightforward:
\begin{align}
\big(T_\mu fg\big)(x)=f\big(T_\mu g\big)(x)+g\big(T_\mu f\big)(x)-4\mu x\left(H_{-1}f\right)(x)\left(H_{-1}g\right)(x).\label{eq1.12}
\end{align}
 Suppose now that $\{P_n\}_{n\geqslant0}$ is  ``$T_\mu$-classical”. Thus $\{Q_n\}_{n\geqslant0}$ is also orthogonal and then satisfies
    a second order recurrence relation.
 Furthermore, if we assume that $\{P_n\}_{n\geqslant0}$ is symmetric, we easily check that $\{Q_n\}_{n\geqslant0}$ is also symmetric, that is,
  \begin{subequations}
\begin{align}
&Q_{n+2}(x)=xQ_{n+1}(x)-{\tilde\gamma}_{n+1}Q_n(x),\ n\geqslant0,\label{eq1.13a}\\
& Q_0(x) = 1,\ Q_1(x) = x.\label{eq1.13b}
\end{align}
    \end{subequations}
  As mentioned above,  it is not so easy to apply our procedure when the PS $\{P_n\}_{n\geqslant0}$ is not symmetric.
 Accordingly we are in effect dealing only with the symmetric case and therefore our starting points are the recurrence relations
 \eqref{eq1.11a}-\eqref{eq1.11b} and \eqref{eq1.13a}-\eqref{eq1.13b}.\par
\noindent First, from \eqref{eq1.11a}-\eqref{eq1.11b} by applying the operator $T_\mu$ we get
   \begin{subequations}
\begin{align}
&\big(T_\mu P_{n+2}\big)(x)=\big(T_\mu xP_{n+1}\big)(x)-\gamma_{n+1}\big(T_\mu P_n\big)(x),\ n\geqslant0,\label{eq1.14a}\\
& \big(T_\mu P_0\big)(x) = 0,\ \big(T_\mu P_1\big)(x) = \mu_1.\label{eq1.14b}
\end{align}
    \end{subequations}
Applying the formula \eqref{eq1.12} for $f(x)=x$ and $g(x)=P_{n+1}(x)$, we see that
 \begin{align}
\big(T_\mu xP_{n+1}\big)(x)=x\big(T_\mu P_{n+1}\big)(x)+\xi_{n+1}P_{n+1},\ n\geqslant0,\label{eq1.15}
\end{align}
where we have put $\xi_{n}=1+2\mu(-1)^{n},\ n\geqslant0$.
Substituting \eqref{eq1.15} into \eqref{eq1.14a} yields
\begin{align}
\xi_{n+1}P_{n+1}(x) = \mu_{n+2}Q_{n+1}(x)+\mu_{n}\gamma_{n+1}Q_{n-1}(x)-\mu_{n+1}xQ_{n}(x),\ n\geqslant0.\label{eq1.16}
\end{align}
Using the $\xi_{n}$-coefficients will be very helpful later on to simplify most operations in the next two sections.
 As one can easily see, the coefficients $\xi_{n}$ and $\mu_{n}$ are interrelated via
\begin{align}
\xi_{n} = \mu_{n+1}-\mu_{n},\ n\geqslant0.\label{eq1.17}
\end{align}
\noindent  At the end  of this section, let us remember the definition of the {\it shifted} polynomials denoted $\{\hat{P}_n\}_{n\geqslant0}$
corresponding to the PS $\{P_n\}_{n\geqslant0}$. For all $n, n=0, 1,\ldots,$ we have
\begin{align}
{\hat P}_n(x) :=a^{-n} P_n({a}x+{b}),\ \hbox{\rm for}\ ({a} ; {b})\in \C^\ast\times\C.\label{eq1.18}
\end{align}
Since the classical character of the considered polynomials is preserved by any linear change of the variable, for the OPS $\{P_n\}_{n\geqslant0}$
satisfying \eqref{eq1.8a}-\eqref{eq1.8b}, we obtain that the polynomials ${\hat P}_n, {n=0, 1,\ldots}$, satisfy also the second-order recurrence relation
 \begin{subequations}
\begin{align}
&{{\hat P}}_{n+2}(x) = (x - {\hat\beta}_{n+1}){{\hat P}}_{n+1}(x)-{\hat\gamma}_{n+1}{{\hat P}}_n (x), \ n\geqslant0,\label{eq1.19a}\\
&{{\hat P}}_1(x) = x - {\hat\beta}_0,\ {{\hat P}}_0(x) = 1,\label{eq1.19b}
\end{align}
 \end{subequations}
with
\begin{align}
{\hat\beta}_{n}=\frac{{\beta}_{n}-{b}}{{a}},\, n\geqslant0,\ \  \mbox{and}\ \ {\hat\gamma}_{n+1}=\frac{{\gamma}_{n+1}}{{a}^2}, \,
 n\geqslant0\quad ({a}\ne0).\label{eq1.20}
\end{align}
Now, under the hypothesis that $\{P_n\}_{n\geqslant0}$ is symmetric, we have ${\hat\beta}_{n}={\beta}_{n}=0,\, n\geqslant0$.\par
\noindent The structure of this paper is as follows. In Section 2, we see how Definition 1.3 allows us to construct the symmetric $T_\mu$-classical
 orthogonal polynomials by solving a system of nonlinear equations satisfied by the coefficients  ${\gamma}_{n}$ and ${\hat\gamma}_{n}$.
 Based on solutions of this system  we state and prove the two  theorems characterizing the symmetric Dunkl-classical orthogonal polynomials.
The third section is devoted to the application of both theorems. Theorem 2.1 asserts that the only symmetric Dunkl-classical orthogonal
 polynomials are the generalized Hermite  polynomials and the generalized Gegenbauer  polynomials, and Theorem 2.2 gives a characterization of these
  families via a differential-difference equation of type \eqref{eq1.2}.
  Many other known properties of each family are given.
  %%%%%%%%%%%%%%%%%%%%%%%%%%%%%
 \section{Construction of the symmetric $T_\mu$-classical orthogonal polynomials}
The main result of the present paper is the following theorem which ensures the existence and the uniqueness  of the symmetric Dunkl-classical orthogonal
polynomial sequences. For this purpose, we show that the coefficients ${\gamma_{n}}$ have a specific rational expression which
actually characterizes the  symmetric Dunkl-classical orthogonal polynomials.
\begin{theorem} Let $\{P_n\}_{n\geqslant0}$ be a MOPS and $\{Q_n\}_{n\geqslant0}$ the sequence defined by $Q_n(x)= \mu_{n+1}^{-1}\big(T_\mu P_{n+1}\big)(x)$.
Then $\{P_n\}_{n\geqslant0}$ satisfies the recurrence relation  \eqref{eq1.11a}-\eqref{eq1.11b} with
\begin{align}
{\gamma_{n+2}}
=\frac{\xi_{n+2}}{\xi_{n+1}}\frac{\mu_{n+2}}{\mu_{n+1}}\frac{\big(\mu_{n}\theta_{n}-\mu_{n-1} \big)}{\big(\mu_{n+3}\theta_{n+1}-
\mu_{n+2}\big)}{\gamma_{n+1}},\ n\geqslant0,\label{eq2.1}
\end{align}
where $\theta_n$ a  solution to the Riccati equation
\begin{align}
\xi_{n+1}\theta_{n+1}+  {\xi_{n}}{\theta_n^{-1}} &=2,\ n\geqslant1,\label{eq2.2}
\end{align}
if and only if $\{Q_n\}_{n\geqslant0}$ is a MOPS satisfying   \eqref{eq1.13a}-\eqref{eq1.13b} with
\begin{align}
\tilde{\gamma}_n=\frac{\mu_{n}}{\mu_{n+1}}\theta_n\gamma_{n+1},\ n\geqslant1.\label{eq2.3}
\end{align}
\end{theorem}
Note that, in the identity \eqref{eq2.1}, we have adopted the  conventions that $\theta_0=1$ and $\mu_{-1}=-1+2\mu$.\\
\begin{proof}
  First we shall prove that if $\{Q_n\}_{n\geqslant0}$ is a MOPS satisfying   \eqref{eq1.13a}-\eqref{eq1.13b} with \eqref{eq2.3}, then
$\{P_n\}_{n\geqslant0}$ satisfies  \eqref{eq1.11a}-\eqref{eq1.11b} with \eqref{eq2.1}-\eqref{eq2.2}.\\
On account of   \eqref{eq1.13a}-\eqref{eq1.13b} with the use of \eqref{eq2.3}, it is easily seen that the equation \eqref{eq1.16} writes
  \begin{subequations}
\begin{align}
P_{n+1}(x) &= Q_{n+1}(x)+\frac{\mu_{n}}{\xi_{n+1}}\gamma_{n+1}\big(1-\theta_n\big)Q_{n-1}(x),\ n\geqslant1,\label{eq2.4a}\\
P_{1}(x) &= Q_{1}(x)=x, \, P_0(x)= Q_0(x)=1. \label{eq2.4b}
\end{align}
  \end{subequations}
  Substituting, in \eqref{eq1.11a}, $P_{n+2}, P_{n+1}$ and $ P_n$ by their expressions provided in \eqref{eq2.4a}, using again \eqref{eq1.13a}
we get an expansion in terms of $Q_n$ and  $Q_{n-2}$. After an easy  computation we derive the equations
\begin{align*}
\frac{\big(\mu_{n+3}\theta_{n+1}-\mu_{n+2}\big)}{\xi_{n+2}}\frac{\gamma_{n+2}}{\mu_{n+2}}
&=\frac{\big(\mu_{n}\theta_{n}-\mu_{n-1} \big)}{\xi_{n+1}}\frac{\gamma_{n+1}}{\mu_{n+1}},\ n\geqslant1,\\
\big(\mu_3\theta_1-\mu_2\big)\gamma_{2}&=\mu_2\gamma_{1};\\
\xi_{n+1}\theta_{n+1}+  {\xi_{n}}{\theta_n^{-1}} &=2,\ n\geqslant1.
\end{align*}
Observe that the last equation  is exactly \eqref{eq2.2}. If moreover we adhere to the above conventions, the first two ones may be grouped together
as in \eqref{eq2.1}. The condition is then sufficient.\\
Conversely, if now $\{P_n\}_{n\geqslant0}$  fulfills   \eqref{eq1.11a}-\eqref{eq1.11b} with  $\gamma_{n+1}$ and $\theta_{n}$ satisfying \eqref{eq2.1}
and \eqref{eq2.2}, we will show that $\{Q_n\}_{n\geqslant0}$ is a symmetrical MOPS. In this part,
 the proof requires tedious computations and its basic idea is similar in spirit to that of \cite[Th. 3.4]{LoVa}.\\
Changing $n$ into $n-1$ in \eqref{eq1.11a}, then substitute $P_{n+1}$, $P_{n}$ and $P_{n-1}$ by their respective expressions provided in \eqref{eq1.16}
 we obtain
\begin{align}
{\xi_{n}} {\mu_{n+2}}Q_{n+1}(x)&=2{\mu_{n+1}}xQ_{n}(x)+2{\mu_{n-1}}\gamma_{n}xQ_{n-2}(x)
-\xi_{n+1}{\mu_{n}}x^2Q_{n-1}(x)\notag\\
&\ -{\xi_{n}}\mu_n\left(\gamma_{n+1}+\gamma_{n}\right)Q_{n-1}(x)-{\xi_{n}}\mu_{n-2}\gamma_{n}\gamma_{n-1}Q_{n-3}(x)\quad (Q_{-n}=0,\, n\geqslant1).
\label{eq2.5}
\end{align}
As the sequence $\{Q_n\}_{n\geqslant0}$ is a basis for the vector space $\mathscr P$, there exist a coefficients $\lambda_{n+1,\nu}$ such that
\begin{align}
xQ_n(x)=\sum_{\nu=0}^{n+1}\lambda_{n+1,\nu}Q_\nu(x), \, n\geqslant0, \label{eq2.6}
\end{align}
where $\lambda_{n+1,n+1}=1$ because $Q_n$ is monic (with $\lambda_{n+1,\nu}=0$ if $\nu<0$ or $\nu>n+1$).\\
 The objective here is to show that $\lambda_{n+1,n-1}\ne0,\,  n\geqslant 1$, and, for all $\nu\ne  n-1$, we have $\lambda_{n+1,\nu}=0$.\\
Use of \eqref{eq2.6} shows that
\begin{align}
x^2Q_n(x)=\sum_{\nu=1}^{n+2}\sum_{\tau=\nu-1}^{n+1}\lambda_{n+1,\tau}\lambda_{\tau+1,\nu}Q_\nu(x)+
\sum_{\tau=0}^{n+1}\lambda_{n+1,\tau}\lambda_{\tau+1,0}Q_0(x), \, n\geqslant0. \label{eq2.7}
\end{align}
 Because of the symmetry property of $\{Q_n\}_{n\geqslant0}$, due to \eqref{eq2.6}, it is easy to check that the coefficients
 $\lambda_{n+1,n-\nu}=0$, for each even number $\nu$,  $\ 0\leqslant \nu\leqslant n$.
  In other words, $\lambda_{n,\nu}=0$ whenever $ n-\nu \equiv 1\ \mbox{\rm mod}\ 2.$\\
 Thus, to prove  \eqref{eq1.13a}-\eqref{eq1.13b}, it only remains to verify that $\lambda_{n+1,n-\nu}=0$, for each odd number $\nu\ne1$.\\
 For this, replace in \eqref{eq2.5} the terms $xQ_{n}$, $xQ_{n-2}$ and $x^2Q_{n-1}$ by their respective expressions obtained from \eqref{eq2.6}
  and \eqref{eq2.7}   to get
\begin{align}
 {\xi_{n}}{\mu_{n+2}}Q_{n+1}(x)&=2{\mu_{n+1}}\sum_{\nu=0}^{n+1}\lambda_{n+1,\nu}Q_\nu(x)+
 2{\mu_{n-1}}\gamma_{n}\sum_{\nu=0}^{n-1}\lambda_{n-1,\nu}Q_\nu(x)\notag\\
&\ \,-\xi_{n+1}{\mu_{n}}\Big[\sum_{\nu=1}^{n+1}\sum_{\tau=\nu-1}^{n}\lambda_{n,\tau}\lambda_{\tau+1,\nu}Q_\nu(x)+
\sum_{\tau=0}^{n}\lambda_{n,\tau}\lambda_{\tau+1,0}Q_0(x) \Big]\notag\\
&\ \,-{\xi_{n}}\mu_n\left(\gamma_{n+1}+\gamma_{n}\right)Q_{n-1}(x)-{\xi_{n}}\mu_{n-2}\gamma_{n}\gamma_{n-1}Q_{n-3}(x).\label{eq2.8}
\end{align}
\noindent From this, for $n\geqslant1$, equating the coefficients of $Q_{n-1}$ leads to
\begin{align}
\xi_{n} \mu_{n+2}\lambda_{n+1,n-1}-\xi_{n+1}\mu_{n} \lambda_{n,n-2} +\xi_{n+1}\mu_{n-2}\gamma_{n}-\xi_{n}\mu_{n}\gamma_{n+1}=0.\label{eq2.9}
\end{align}
 In particular, if we replace $n=1$ in  the last equation we find
\begin{align*}
  \mu_3 \lambda_{2,0}-  \mu_{1}\left(\gamma_{2}+\gamma_{1}\right)=0.
\end{align*}
But, for $n=0$ in \eqref{eq2.1}, we can write $\gamma_{1}=(\mu_{3}\theta_{1}-2) \frac{\gamma_{2}}{\mu_{2}}$ giving
\begin{align*}
  \lambda_{2,0}= \frac{\mu_{1}}{\mu_{2}}\theta_{1} \gamma_{2}.
\end{align*}
We now proceed by induction and assume   that the following formula holds
\begin{align}
\lambda_{n,n-2}=\frac{\mu_{n-1}}{\mu_{n}}\theta_{n-1} \gamma_{n}, \ \ \mbox{ for}\  n\geqslant2.\label{eq2.10}
\end{align}
Because of \eqref{eq2.1} and \eqref{eq2.10},   we easily deduce from \eqref{eq2.9} that
\begin{align*}
 \lambda_{n+1,n-1} =\frac{\mu_{n}}{\mu_{n+1}} \theta_{n}\gamma_{n+1}\ne0, \ n\geqslant 1.
\end{align*}
Similarly,  for $n\geqslant3$, equating the coefficients of $Q_{n-3}$ gives
\begin{align*}
\xi_{n} \mu_{n+2}\lambda_{n+1,n-3}-\xi_{n+1}\mu_{n}\lambda_{n,n-4}+\left(2\mu_{n-1}\gamma_{n}
-\xi_{n+1}\mu_{n}\lambda_{n,n-2}\right)\lambda_{n-1,n-3} -\xi_{n}\mu_{n-1}\gamma_{n}\gamma_{n-1}=0.
\end{align*}
After replacing  $\lambda_{n,n-2}$ and $\lambda_{n-1,n-3}$ by their expressions provided by \eqref{eq2.10}, the latter becomes
\begin{align}
\xi_{n} \mu_{n+2}\lambda_{n+1,n-3}-\xi_{n+1}\mu_{n}\lambda_{n,n-4} =
\mu_{n-2}\gamma_{n}\gamma_{n-1}\big[\xi_{n+1}\theta_{n-1}\theta_{n-2}-2\theta_{n-2}+\xi_{n}\big].\label{eq2.11}
\end{align}
Since $\xi_{n+1}=\xi_{n-1}$  and so  $\xi_{n}=\xi_{n-2}$,  the  Riccati equation \eqref{eq2.2} readily implies the right-hand side of
the latter equality to be zero and therefore we find
\begin{align}
\xi_{n} \mu_{n+2}\lambda_{n+1,n-3}=\xi_{n+1}\mu_{n}\lambda_{n,n-4}, \, n\geqslant3.\label{eq2.12}
\end{align}
If we take $n=3$ in \eqref{eq2.12}, we get $\lambda_{4,0}=0$, and recursively we deduce that $\lambda_{n+1,n-3}=0$, for all  $n\geqslant3$, or,
equivalently, $\lambda_{n+4,n}=0$, for all  $n\geqslant0$.\\
We can now proceed analogously to show that $\lambda_{n+1,n+1-2k}=0$, for each  $k\geqslant3$.\\
Equating the coefficients of $Q_{n+1-2k}$ for $k\geqslant3$ in \eqref{eq2.8} leads to
\begin{align*}
\xi_{n}\mu_{n+2}\lambda_{n+1,n+1-2k}+\left(2\mu_{n-1}\gamma_{n}-\xi_{n+1}\mu_{n}\lambda_{n,n-2}\right)\lambda_{n-1,n+1-2k}
-\xi_{n+1}\mu_{n}\!\!\sum_{\tau=n-2k}^{n-6}\lambda_{n,\tau}\lambda_{\tau+1,n+1-2k}=0.
\end{align*}
Based on the identity \eqref{eq2.10}, this may be  also written as
\begin{align*}
\xi_{n}\mu_{n+2}\lambda_{n+1,n+1-2k}+\mu_{n-1}\gamma_{n}\left(2-\xi_{n+1}\theta_{n-1}\right)\lambda_{n-1,n+1-2k}
-\xi_{n+1}\mu_{n}\sum_{\tau=3}^{k}\lambda_{n,n-2\tau}\lambda_{n+1-2\tau,n+1-2k}=0.
\end{align*}
Suppose that $\lambda_{n,n-2j}=0$ for each $j=2, \ldots, k-1$, $n\geqslant 2j$. This readily gives
\begin{align}
\xi_{n}\mu_{n+2}\lambda_{n+1,n+1-2k}=\xi_{n+1}\mu_{n}\lambda_{n,n-2k}.\label{eq2.13}%\label{eq1.30}
\end{align}
From the above, by decreasing $n$, it follows that
\begin{align}
 \lambda_{n+1,n+1-2k}=\frac{\xi_{n+1}\mu_{2k+1}\mu_{2k}}{\xi_{2k}\mu_{n+2}\mu_{n+1}}\lambda_{2k,0}.\label{eq2.14}
\end{align}
In particular, for $n=2k$, we have
\begin{align}
 \lambda_{2k+1,1}=\frac{\xi_{2k+1}\mu_{2k}}{\xi_{2k}\mu_{2k+2}}\lambda_{2k,0}.\label{eq2.15}%\label{eq1.32}
\end{align}
Finally, for $n\geqslant6$, a simple comparison of the coefficients of $Q_0$ in \eqref{eq2.8} leads to
\begin{align}
\xi_{n}\mu_{n+1}\lambda_{n+1,0}+2\mu_{n-1}\gamma_{n}\lambda_{n-1,0}-\xi_{n+1}\mu_{n}\sum_{\tau=0}^{n-1}\lambda_{n,\tau}\lambda_{\tau+1,0}=0.%
\label{eq2.16}%\label{eq1.33}
\end{align}
But, due to the symmetry property of $\{Q_n\}_{n\geqslant0}$, we have $\lambda_{2n+1,0}=0$ for $n\geqslant0$. Accordingly,
if we replace $n$ in \eqref{eq2.16} by $2n+1$, this identity may be written as
\begin{align*}
\xi_{2n+1}\mu_{2n+2}\lambda_{2n+2,0}+2\mu_{2n}\gamma_{2n+1}\lambda_{2n,0}-\xi_{2n+2}\mu_{2n+1}\lambda_{2n+1,1}\lambda_{2,0}
-\xi_{2n+2}\mu_{2n+1}\sum_{\tau=1}^{n-2}\lambda_{2n+1,2\tau+1}\lambda_{2\tau+2,0}=0.
\end{align*}
Based on the identities \eqref{eq2.14}-\eqref{eq2.15}, if we suppose that $\lambda_{2k,0}=0$, $k=2, \ldots, n$, we deduce by induction that
$\lambda_{2n+2,0}=0$, and then obtain
\begin{align*}
 \lambda_{n+1,n+1-2k}=0, \quad \mbox{for any }\  k\geqslant 3.
\end{align*}
As a consequence of all the above results, we deduce that \eqref{eq2.6} writes
\begin{align*}
xQ_n(x)=Q_{n+1}+\lambda_{n+1,n-1}Q_{n-1}(x), \, n\geqslant0,  \quad \mbox{with }\  \lambda_{n+1,n-1} =\frac{\mu_{n}}{\mu_{n+1}} \theta_{n}\gamma_{n+1}\ne0.
\end{align*}
This means that $\{Q_n\}_{n\geqslant0}$ is a symmetrical OPS which is the desired conclusion.
\end{proof}
{\bf Remark 1.}
Clearly, the identification of the polynomial sequences discussed here depends on the solutions of the Riccati equation \eqref{eq2.2}, which are given by
\begin{align*}
  {\bf A.} \quad\theta_n &= 1,\ \ n\geqslant1, \quad\mbox{(the trivial solution)};\\
 {\bf B.}\quad\theta_n &= \frac{n + \theta + 1- \mu(-1)^{n}}{n +\theta +\mu(-1)^{n}},\ \ n\geqslant1,\quad \mbox{with} \ \theta \ne - n\pm\mu(-1)^{n}
 \ \ \mbox{an   arbitrary parameter}. \hskip5cm
\end{align*}
In consequence, without imposing  any additional conditions, we can then say that there are exactly two families of symmetric Dunkl-classical
OPS that are solutions to our problem: the generalized Hermite polynomials (obtained in Case {\bf A}) and the generalized Gegenbauer
polynomials (obtained in Case {\bf B}).\\
 It is worth mentioning that when $\mu=0$, the Dunkl operator coincides with the usual differential operator, say $T_0 = D$, and so
 the system \eqref{eq2.1}-\eqref{eq2.2} is  simplified to
 \begin{align*}
&{\gamma_{n+2}}
=\frac{{n+2}}{{n+1}}\frac{\big(n(\theta_{n}-1)+1 \big)}{\big({(n+3)}(\theta_{n+1}-1)+1\big)}{\gamma_{n+1}},\ n\geqslant0,\\
&\theta_{n+1}+  {\theta_n^{-1}} =2,\ n\geqslant1.
\end{align*}
This system was already solved providing the two  families of symmetric  classical orthogonal polynomials: Hermite polynomials and
Gegenbauer polynomials (see \cite{BeDo1} for more details).\\
Before discussing the corresponding results when $\mu\ne0$ with the identification  of the obtained OPS as announced above,
 we first give the next theorem which is an important consequence of Theorem 2.1.
 Note that a partial result of  this theorem was presented in \cite[Lemma 2.1]{BeGa}, however, our proof is
quite different from that of this paper. Here, we also prove  the reciprocal condition found in this work.
\begin{theorem} Let $\{P_n\}_{n\geqslant0}$ be a symmetric MOPS satisfying the recurrence relation  \eqref{eq1.11a}-\eqref{eq1.11b}.
Then $\{P_n\}_{n\geqslant0}$ is $T_\mu$-classical if and only if each polynomial $P_n$ satisfies the differential-difference equation
\begin{align}
\left(a_n x^2-b_n\right)T^2_\mu P_n-c_nxT_\mu P_n-d_nP_n=0,\label{eq2.17}
\end{align}
where
  \begin{subequations}
\begin{align}
a_n&= \theta_{n-1}-1,\label{eq2.18a}\\
b_n&= {\left(\mu_{n+1}\theta_{n-1}-\mu_{n}\right)\left(\mu_{n-1}\theta_{n-1}-\mu_{n-2}\right)}\frac{\gamma_{n}}{\xi_{n}\mu_{n}},\label{eq2.18b}\\
c_n&= \mu_{n-2}\theta_{n-1}-\mu_{n-1},\label{eq2.18c}\\
d_n&= \xi_{n}\mu_{n}\theta_{n-1},\label{eq2.18d}
\end{align}
  \end{subequations}
  with $\theta_n$ a solution to  the Riccati equation \eqref{eq2.2}.
  It is understood that $\theta_{-1}=0$.
\end{theorem}
Before proving Theorem 2.2, it must be noted that the expression defining $b_n$ is initially given by
\begin{equation}
b_n= {\left(\mu_{n+1}\theta_{n-1}-\mu_{n}\right)\left(\mu_{n+2}\theta_{n}-\mu_{n+1}\right)}\frac{\gamma_{n+1}}{\xi_{n+1}\mu_{n+1}}.\label{eq2.19}
\end{equation}
Indeed, after the substitution of the coefficient $\gamma_{n+1}$ by an expression  in terms of $\gamma_{n}$
 provided in \eqref{eq2.1} we see at once that  \eqref{eq2.19} is an equivalent form of \eqref{eq2.18b}.\\
\begin{proof}
  The  necessary condition. If we replace $n$ in \eqref{eq1.13a} by $n-1$, taking into account \eqref{eq2.3}, we get
   \begin{align}
\mu_{n+1}Q_{n+1}(x)=\mu_{n+1}xQ_{n}(x)-\mu_{n}\theta_n\gamma_{n+1}Q_{n-1}(x),\ n\geqslant0,\label{eq2.20}
\end{align}
   On combining this with \eqref{eq1.16} yields
\begin{align}
\xi_{n+1}\mu_{n+1}P_{n+1}(x)+\mu_{n}\left(\mu_{n+2}\theta_{n}-\mu_{n+1}\right)\gamma_{n+1}Q_{n-1}(x)- \xi_{n+1}\mu_{n+1}xQ_{n}(x)=0,\ n\geqslant0.
\label{eq2.21}
\end{align}
Multiply  \eqref{eq2.21} by $\theta_n$ and use \eqref{eq2.20}, then  shift $n\rightarrow n-1$, we obtain
\begin{align}
\xi_{n}\theta_{n-1}P_{n}(x)+\mu_{n}\left(\theta_{n-1}-1\right)xQ_{n-1}(x)-\left(\mu_{n+1}\theta_{n-1}-\mu_{n}\right)Q_{n}(x)=0,\ n\geqslant0.\label{eq2.22}
\end{align}
If we multiply now \eqref{eq2.21} by $\left(\mu_{n+1}\theta_{n-1}-\mu_{n}\right)$ and \eqref{eq2.22} by $\xi_{n+1}\mu_{n+1}x$ and add the resulting identities,
 it follows that
\begin{align*}
\left(\mu_{n+1}\theta_{n-1}-\mu_{n}\right)P_{n+1}(x)&=\xi_{n}\theta_{n-1}x P_{n}(x)+
 \left(\theta_{n-1}-1\right)x^2T_\mu P_{n}(x)\notag\\
 &-\left(\mu_{n+2}\theta_{n}-\mu_{n+1}\right)\left(\mu_{n+1}\theta_{n-1}-\mu_{n}\right)\frac{\gamma_{n+1}}{\xi_{n+1}\mu_{n+1}}T_\mu P_{n}(x).%\label{eq2.22}
\end{align*}
The action of $T_\mu$ on both sides of the latter gives us
\begin{align}
\left(\mu_{n+1}\theta_{n-1}-\mu_{n}\right)T_\mu P_{n+1}(x)&=\left(a_n x^2-b_n\right)T^2_\mu P_n+\left((\xi_{n}+2)\theta_{n-1}-2\right)
xT_\mu P_n+\xi_{n}^2\theta_{n-1} P_{n}(x).\label{eq2.23}
\end{align}
On the other hand, we see that \eqref{eq2.22} writes also in the form
\begin{align}
\left(\mu_{n+1}\theta_{n-1}-\mu_{n}\right)T_\mu P_{n+1}(x)=
\mu_{n+1}\xi_{n}\theta_{n-1}P_{n}(x)+\mu_{n+1}\left(\theta_{n-1}-1\right)xT_\mu P_{n}(x).\label{eq2.24}
\end{align}
Finally,  equating the right hand sides of both \eqref{eq2.23}  and \eqref{eq2.24}, the differential-difference equation \eqref{eq2.17} follows
 immediately with the coefficients $a_n,\ldots, d_n$ are defined by \eqref{eq2.18a}-\eqref{eq2.18d}.\\
Conversely, suppose that the MOPS $\{P_n\}_{n\geqslant0}$ fulfils \eqref{eq2.17} with the definitions \eqref{eq2.18a}-\eqref{eq2.18d}.
Our  objective is to establish the orthogonality of the sequence $\{Q_n\}_{n\geqslant0}$.\\
Since the polynomials  $P_n$ are monic and symmetric, if we write $P_n(x)=x^n+\alpha_{n}x^{n-2}+\cdots$, for  $n\geqslant2$, and insert this in
 \eqref{eq1.11a}, make the change $n\rightarrow n-1$, a simple identification gives rise to
\begin{align}
\gamma_{n}=\alpha_{n}-\alpha_{n+1},\ n\geqslant1.\label{eq2.25}
\end{align}
Besides, substitute the above expansion of $P_n(x)$  in the differential-difference equation \eqref{eq2.17}, we easily check that the coefficients
 of $x^n$ and $x^{n-2}$ provide the following equalities
\begin{align}
d_{n}=\mu_{n}\mu_{n-1}a_n-\mu_{n}c_n,\ n\geqslant0,\label{eq2.26}
\end{align}
and
\begin{align}
\big(\mu_{n-2}\mu_{n-3}a_n-\mu_{n-2}c_n-d_n\big)\alpha_{n}=\mu_{n}\mu_{n-1}b_n,\ n\geqslant0.\label{eq2.27}
\end{align}
\noindent On combining \eqref{eq2.26}  with \eqref{eq2.27}, taking into consideration \eqref{eq2.18a}-\eqref{eq2.18d}, we see that
\begin{align}
\alpha_{n}=-\frac{\mu_{n}\mu_{n-1}b_n}{2\big(\mu_{n-1}\theta_{n-1}-\mu_{n-2}\big)},\ n\geqslant0.\label{eq2.28}
\end{align}
Therefore \eqref{eq2.25} is readily written as
\begin{align}
\gamma_{n}=\frac{\mu_{n+1}\mu_{n}b_{n+1}}{2\big(\mu_{n}\theta_{n}-\mu_{n-1}\big)}-
\frac{\mu_{n}\mu_{n-1}b_n}{2\big(\mu_{n-1}\theta_{n-1}-\mu_{n-2}\big)},\ n\geqslant1.\label{eq2.29}
\end{align}
Now, we  substitute in the latter  $b_{n+1}$ (resp. $b_n$) by its respective expression provided from the identity \eqref{eq2.18b}
 (resp.\eqref{eq2.19}) obtaining
\begin{align}
2\xi_{n+1}\mu_{n+1}\big(\mu_{n-1}\theta_{n-1}-\mu_{n-2}\big)\gamma_{n}
=&\left[\mu_{n+1}\big(\mu_{n-1}\theta_{n-1}-\mu_{n-2}\big)-\mu_{n-1}\big(\mu_{n+1}\theta_{n-1}-\mu_{n}\big)\right]\times\notag\\
&\ \mu_{n}{\big(\mu_{n+2}\theta_{n}-\mu_{n+1}\big)}\gamma_{n+1},\ n\geqslant1.\label{eq2.30}
\end{align}
An easy computation shows that the expression between square brackets simplifies as
\[
\mu_{n+1}\big(\mu_{n-1}\theta_{n-1}-\mu_{n-2}\big)-\mu_{n-1}\big(\mu_{n+1}\theta_{n-1}-\mu_{n}\big)=2\xi_{n}.
\]
Substituting this into \eqref{eq2.30} and then replacing $n$ by $n+1$ we conclude that
\begin{align*}
{\gamma_{n+2}}
=\frac{\xi_{n+2}}{\xi_{n+1}}\frac{\mu_{n+2}}{\mu_{n+1}}\frac{\big(\mu_{n}\theta_{n}-\mu_{n-1} \big)}{\big(\mu_{n+3}\theta_{n+1}-
\mu_{n+2}\big)}{\gamma_{n+1}},\ n\geqslant0,
\end{align*}
and so we have found the equality \eqref{eq2.1}. In consequence, use of Theorem 2.1 ensures the orthogonality of the sequence $\{Q_n(x)\}$,
 which means that $\{P_n(x)\}$ is $T_\mu$-classical, and the proof is complete.
\end{proof}
%%%%%%%%%%%%%%%%%%%%%%%%%%%%%%%%%%%%
{\bf Remark 2.} As we  observed in Remark 1,  for $\mu=0$, we again meet  each of the families of Hermite polynomials and Gegenbauer polynomials.
This will be  more  specified  in the next section when applying Theorem 2.2. Equation \eqref{eq2.17} is then reduced  to the ordinary differential
equation of type
\begin{align}
a(x) y''+b(x)y'-\lambda_ny=0,\ \ y=P_n(x),\label{eq2.31}
\end{align}
where $a(x)$ is a $0$-degree polynomial or a $2$-degree polynomial (depending on the cases {\bf A} or {\bf B}), $b(x)$ is a polynomial with degree $1$
 and $\lambda_n$ are scalars.
 %%%%%%%%%%%%%%%%%%%%%%%%%%%%%%%%%%%
\section{Applications }
We illustrate Theorems 2.1 and 2.2 by considering separately the two cases {\bf A} and {\bf B} as referred to above. We confirm that
the only symmetric Dunkl-classical OPS are the generalized Hermite polynomials and the generalized Gegenbauer polynomials, as one may have expected.
 As we have seen, this result was obtained directly from the solutions of Equation \eqref{eq2.2} without imposing any further restrictions.\\
 \noindent In addition to illustrating how our theorems provide the above results, we also remind  some of the known properties satisfied by
 the resulting polynomials. We use, for the most part, the notations of \cite{Chih2}.\par
\medskip
\noindent{\bf Case A.} Application of Theorem 2.1 when $\theta_n = 1$ shows that Equality \eqref{eq2.1} is reduced to
\begin{align*}
\gamma_{n+2}=\frac{\mu_{n+2}}{\mu_{n+1}}\gamma_{n+1},\ n\geqslant0,
\end{align*}
from which, taking into account \eqref{eq2.3}, we immediately deduce that
$\tilde{\gamma}_{n+1}=\gamma_{n+1}=\frac{\gamma_{1}}{\mu_1}{\mu_{n+1}},\ n\geqslant0$.\\
Letting $\delta_n=\mu(1-(-1)^n)$ and choosing $\gamma_{1}=\frac{1}{2}\mu_1$, it follows that
\begin{align}
\tilde{\gamma}_{n+1}=\gamma_{n+1}=\frac{1}{2}\big(n+1+\delta_{n+1}\big),\ n\geqslant0.\label{eq3.1}
\end{align}
From this we see that the regularity conditions still hold provided $\mu\ne-n-\frac{1}{2}, n=0, 1, 2, \ldots$.
Moreover, it follows from \eqref{eq3.1} that $Q_n=P_n,\, n\geqslant0$, which reads as $T_\mu P_{n+1}=\mu_{n+1}P_n,\, n\geqslant0$.
Whence the OPS $\{P_n\}_{n\geqslant0}$ belongs to the class of Appell sequences relative to $T_\mu$. It said to be a Dunkl-Appell sequence.\\
We thus encounter the generalized Hermite polynomials $\{H^{\mu}_n\}_{n\geqslant0}$ which have been studied extensively by Chihara in his Ph.D. thesis
\cite{Chih1}. The same author  was also established in \cite{Chih2} that these polynomials are related to the generalized Laguerre polynomials
$L_n^{\alpha}(x),\, n=0, 1, \ldots$, via
 \begin{align*}
H^{\mu}_{2n}(x) &=  (-1)^{n}2^{2n}n!L_n^{\alpha}(x^{2}),\ n\geqslant0,\\
H^{\mu}_{2n+1}(x) &=  (-1)^{n}2^{2n+1}n!xL_n^{\alpha+1}(x^{2}),\ n\geqslant0,\quad \alpha=\mu-\frac{1}{2}.
\end{align*}
In particular, for $\mu=0$, we have $H^{0}_{n}(x) =  {H}_{n}(x), n\geqslant0$, the classical Hermite polynomials.\\
 The normalization coefficient is chosen so that the leading coefficient of $H^{\mu}_{n}(x)$ becomes $2^n$ and so $H^{\mu}_{n}(x)=2^n \hat{H}^{\mu}_{n}(x)$,
which means that the sequence $\{P_n\}_{n\geqslant0}$ coincides with $\{\hat{H}^{\mu}_n\}_{n\geqslant0}$. \\
 Chihara pointed out that these polynomials were first introduced by Szeg\"{o} in \cite[p.380]{Szeg}  who gave the associated orthogonalizing weight
 $w(x)=|x|^{2\mu}e^{-x^2}$, $\mu>-\frac{1}{2}$, supported on $(-\infty,+\infty)$ and found that each polynomial satisfies two second order differential
 equations. For instance, the following
 \begin{align}
xy''+2(\mu-x^2)y'+(2nx-\delta_nx^{-1})y=0. \label{eq3.2}
\end{align}
  Many other formulas may be found in \cite[pp.156-158]{Chih2}. On the other hand, Maroni \cite{Maro1} obtained that the generalized Hermite polynomials
  satisfy the structure relation
 \begin{align*}xP'_{n+1}(x)=-\delta_{n+1}P_{n+1}(x)+2\gamma_{n+1}xP_{n}(x), \quad P_n(x)=\hat{H}^{\mu}_n(x),\end{align*}
which can also be turned into
\begin{align*}
xP'_{n+1}(x)=(n+1)P_{n+1}(x)+\frac{1}{2}\left(n^2+(2\mu+1)n+\delta_n\right)\!P_{n-1}(x).
\end{align*}
For $\mu=0$, one sees that $Q_n(x)=P_{n}(x)=(n+1)^{-1}P'_{n+1}(x)$. In consequence,  the latter becomes
 \begin{align*}
xQ_{n}(x)=Q_{n+1}(x)+\frac{1}{2}nQ_{n-1}(x),\ n\geqslant0,
\end{align*}
which is clearly the recurrence relation satisfied by the derivatives of the classical Hermite polynomials.\par
\noindent Among the various works concerning applications of generalized Hermite polynomials, we quote, for instance, \cite{Rose,Rosl,ShChF}
(see also the references given there).\par
\medskip
\noindent Now let $\mathscr H_\mu$ be  the linear functional associated to the  polynomials ${H}^{\mu}_{n}, n\geqslant0$.
The regularity of this functional is guaranteed  provided
$\mu\ne-\frac{1}{2}, -\frac{3}{2}, -\frac{5}{2}, \ldots$, and it is  positive definite only when $\mu>-\frac{1}{2}$.\\
By considering some special cases of semi-classical forms, it was  observed in \cite{Maro1} that $\mathscr H_\mu$ is actually a semi-classical
 functional of class one satisfying the distributional differential  equation \eqref{eq1.9},
where $\phi(x)=x$ and $\psi(x)=2x^2-(2\mu+1)$. This clearly  specifies that $s=\max\left(\deg\phi-2,\deg\psi-1\right)=1$.\\
As above-mentioned, for $\mu>-\frac{1}{2}$,  this functional is  represented as
\begin{align*}
\big<\mathscr H_\mu,f\big>=\int_{-\infty}^{+\infty}f(x)|x|^{2\mu}e^{-x^2}dx, \ \ \forall f \in \mathscr P,
\end{align*}
so that its moments are given by
 \begin{subequations}
 \begin{align}
&(\mathscr H_\mu)_{2n}={\Gamma\big(n+\mu+{1}/{2}\big)}, \ n\geqslant0,\label{3.3a}\\
&(\mathscr H_\mu)_{2n+1}=0, \ n\geqslant0.\label{3.3b}
\end{align}
 \end{subequations}
When $\mu<-\frac{1}{2}$, however, finding orthogonalizing weight function is more difficult. Morton and Krall \cite{MoKr} attacked this problem from
 the point of view of distribution theory. In fact, based on the canonical regularization distributional technique introduced in \cite[Ch.1, Sec.3]{GeSh},
 the authors had developed a procedure for constructing distributional weight functions supported on  a succession of disjointed strips.\\
Applying this method, Krall showed in \cite{Kral}  that an integral representation of the functional $\mathscr H_\mu$  can be obtained,
 in any strip of the form $-j-\frac{1}{2}<\mu<-j+\frac{1}{2}$,  where $j$ is a fixed positive integer, in the form
  \begin{align}
\big<\mathscr H_\mu,f\big>=\int_{0}^{+\infty}\!\!|x|^{2\mu}\Big\{e^{-x^2}[f(x)+f(-x)]-
\!\!\sum_{k=0}^{2j-2}\bigl(e^{-x^2}[f(x)+f(-x)]\bigr)^{(k)}(0)\frac{x^{k}}{k!}\Big\}dx,
\label{eq3.4}
\end{align}
 where $f$ is a sufficiently differentiable function and $g^{(k)}$ stands for the $k^{th}$ derivative of the function $g$.
Note that in this case there is a double regularization of $|x|^{2\mu}$ at $x = 0$ following Gelfand and Shilov \cite{GeSh}.\\
The $2j$-fold integration by parts in \eqref{eq3.4} shows that
\begin{align}
\big<\mathscr H_\mu,f\big>=\frac{\Gamma(2\mu+1)}{\Gamma(2\mu+2j+1)}
\int_{-\infty}^{+\infty}\!|x|^{2\mu+2j}\big(e^{-x^2}f(x)\big)^{(2j)}dx,\ \ j=1, 2, \ldots.\label{eq3.5}
\end{align}
A third form of the distributional weight function involving the classical Hermite polynomials has also been  given in \cite{Kral}.
 This representation is fully  equivalent to both \eqref{eq3.4} and \eqref{eq3.5} when $f$ is a polynomial.
In doing so, the author had taken as the only starting point the moments \eqref{3.3a}-\eqref{3.3b}, under the regularity conditions
imposed on  $\mu$. For more details about this procedure we refer the reader to \cite{MoKr,Kral}.\par
\medskip
\noindent  We continue in this fashion with  Theorem 2.2 to obtain that \eqref{eq2.18a}-\eqref{eq2.18b} simplify to
$$
a_n=0, \, b_n=\xi_n/2, \, c_n=-\xi_n\ \, \mbox{and}\ \, d_n=\xi_n\mu_n.
$$
 On substituting these into \eqref{eq2.17}  we obtain
\begin{align}
T^2_\mu P_n-2xT_\mu P_n+\lambda_nP_n=0,\label{eq3.6}
\end{align}
where $\lambda_n=2\mu_n=2\big(n+\delta_n\big), \ n\geqslant0, $ or, equivalently,
\begin{align*}
\lambda_{2n}=4n, \ n\geqslant0,\ \ \mbox{and} \ \
\lambda_{2n+1}=2\big(2n+2\mu+1\big), \ n\geqslant0.
\end{align*}
For $\mu=0$, we easily check that  \eqref{eq3.6}, \eqref{eq3.2} are identical and coincide with \eqref{eq2.31}, where
$a(x)=1$, $b(x)=-2x$ and $\lambda_n=-2n$. It follows that the equation outlined here is the second-order differential equation
satisfied by  the ordinary Hermite polynomials $H_{n}, \, n=0, 1, \ldots$.\\
 \bigskip
\noindent{\bf Case B.} Similarly to the previous case, on account of the  value  of $\theta_n$, we can write \eqref{eq2.1} in the form
 \begin{align}
{\gamma_{n+2}}=\frac{\big(2n+ \theta+\mu\big)\big(n+1+\theta-\mu(-1)^{n}\big)}
{\big(2(n+2)+\theta+\mu\big)\big(n+\theta+\mu(-1)^{n}\big)}\frac{\mu_{n+2}}{\mu_{n+1}}{\gamma_{n+1}},\ n\geqslant0.\label{eq3.7}
\end{align}
   Multiply, in \eqref{eq3.7},  the numerator and denominator of the right-hand side by $2(n+1)+\theta+\mu$ to readily deduce by telescopy  that
  \begin{align*}
\gamma_{n+1}=\big(\theta+\mu+2\big)\frac{\gamma_1}{\mu_1}
\frac{\big(n+ \theta+\mu(-1)^{n}\big)\big(n+1+\delta_{n+1}\big)}{\big(2n+2+\theta+\mu\big)\big(2n+\theta+\mu\big)},\ n\geqslant0.
\end{align*}
For two fixed parameters $\alpha$ and $\beta$, if we let $\gamma_1=\frac{\beta+1}{\alpha+\beta+2}$ and set
$\, \theta+\mu=2\big(\alpha+\beta+1\big),\ \theta-\mu=2\alpha+1$,  we get  $2\mu=2\beta +1$ and $\mu_1=2(\beta+1)$.
Thus, the above expression may be rewritten  in the form \cite{Belm}:
   \begin{align}
\gamma_{n+1}=
\frac{\big(n+1+2\alpha+\delta_{n+1}\big)\big(n+ 1+\delta_{n+1}\big)}{4\big(n+2+\alpha+\beta\big)\big(n+1+\alpha+\beta\big)},\ n\geqslant0.\label{eq3.8}
\end{align}
With the above reparameterization,  we must take into account that the validity of  \eqref{eq3.8} presumes the parameters $\alpha$ and $\beta$
to take values in $\C$ within the range of regularity.\\
The polynomials arising in this case are clearly the generalized Gegenbauer ones which we denote here by $\{{G}^{(\alpha,\beta)}_n\}_{n\geqslant0}$.
 These polynomials were first investigated by Chihara \cite{Chih2} who obtained their main properties by employing  a direct relation between them and
  the Jacobi OP $P_n^{(\alpha,\beta)},\ n\geqslant0$:
 \begin{align*}
{G}^{(\alpha,\beta)}_{2n}(x) &=  P_n^{(\alpha,\beta)}(2x^{2}-1),\ n\geqslant0,\\
{G}^{(\alpha,\beta)}_{2n+1}(x) &= xP_n^{(\alpha,\beta+1)}(2x^{2}-1),\ n\geqslant0.
\end{align*}
With the notations $\hat{P}_n^{(\alpha,\beta)}$ and $\hat{G}_n^{(\alpha,\beta)}$ for the monic Jacobi polynomials and generalized Gegenbauer polynomials, respectively, we have
  \begin{align}
\hat{G}^{(\alpha,\beta)}_{2n}(x) =  {2^{-n}}\hat{P}_n^{(\alpha,\beta)}(2x^{2}-1)\quad \mbox{and}\ \
\hat{G}^{(\alpha,\beta)}_{2n+1}(x) =   {2^{-n}}\hat{P}_n^{(\alpha,\beta+1)}(2x^{2}-1). \label{eq3.9}
\end{align}
So, the resulting polynomials $\{P_n\}_{n\geqslant0}$ obtained in this second case coincide with $\{\hat{G}^{(\alpha,\beta)}_n\}_{n\geqslant0}$,
  and  are therefore orthogonal with respect to the weight function $w(x)=|x|^{2\mu}(1-x^2)^\alpha$  on $(-1,+1)$.\\
Furthermore,  combining \eqref{eq2.3} with \eqref{eq3.8}, we see at once  that the coefficients $\tilde{\gamma}_{n}, \, n\geqslant1$,
defining the sequence $\{Q_n\}_{n\geqslant0}$ are given by
 \begin{align}
 \tilde{\gamma}_{n+1}=
\frac{\big(n+3+2\alpha+\delta_{n+1}\big)\big(n+ 1+\delta_{n+1}\big)}{4\big(n+3+\alpha+\beta\big)\big(n+2+\alpha+\beta\big)},\ n\geqslant0.\label{eq3.10}
\end{align}
Now, a simple comparison between \eqref{eq3.8} and \eqref{eq3.10}  shows that
 \begin{align}
\tilde{\gamma}_{n+1}(\alpha,\beta)=\gamma_{n+1}(\alpha+1,\beta),\ n\geqslant0,\label{eq3.11}
\end{align}
where $\gamma_{n+1}:=\gamma_{n+1}(\alpha,\beta)$ and $\tilde{\gamma}_{n+1}:=\tilde{\gamma}_{n+1}(\alpha,\beta)$. So, with the notations
$P_n(x):=P_n(x;\alpha,\beta),\, n\geqslant0$, and $Q_n(x):=Q_n(x;\alpha,\beta),\, n\geqslant0$,
 we have
$Q_n(x;\alpha,\beta)=P_n(x;\alpha+1,\beta),\, n\geqslant0$, which means that
   \begin{align}
T_\mu P_{n+1}(x;\alpha,\beta)=\mu_{n+1}P_n(x;\alpha+1,\beta),\, n\geqslant0. \label{eq3.12}
\end{align}
A number of works  in different contexts  have been devoted to the  generalized Gegenbauer polynomials (see for instance \cite{Kono,ChIs,Dett,Belm})
highlighting several properties of these polynomials.
For example, it was obtained in  \cite{Belm}  that these polynomials  satisfy the structure relation
 \begin{align*}
x(x^2-1)P'_{n+1}(x)=\left((n+1)x^2+\delta_{n+1}\right)P_{n+1}(x)
-2\gamma_{n+1}(n+\alpha+\beta+2)xP_{n}(x),
\end{align*}
where the coefficients $\gamma_n,\, n\geqslant1$, are given by \eqref{eq3.8}. An equivalent formulation of this relation  is
 \begin{align*}x(x^2-1)P'_{n+1}(x)&=(n+1)P_{n+3}(x)-\left[(n+1)\gamma_{n+2}+\big(3(n+1)+2(\alpha+\beta+1)\big)\gamma_{n+1}-\delta_{n+1}\right]P_{n+1}(x)\\
&\,+\gamma_{n+1}\gamma_{n}\big(3(n+1)+2(\alpha+\beta+1)\big)P_{n-1}(x).
\end{align*}
The same author  also showed that each polynomial satisfies the second order differential equation
 \begin{align}
 x^2(x^2-1)y''+x\left[(2\alpha+2\beta+3)x^2-2\beta-1\right]y'-\left[n(n+2\alpha+2\beta+2)x^2-\delta_{n}\right]y=0.\label{eq3.13}
 \end{align}
If moreover $\mathscr G_{\alpha,\beta}$ denotes the Gegenbauer linear functional associated to the generalized Gegenbauer polynomials,
 he found that $\mathscr G_{\alpha,\beta}$ satisfies \eqref{eq1.9},
where $\phi(x)=x(x^2-1)$ and $\psi(x)=-\left(2(\alpha+\beta+2)x^2+\beta+1\right)$, with $\alpha,\beta>-1$ and $\beta\ne-\frac{1}{2}$.
From this, we see that $\mathscr G_{\alpha,\beta}$ is  also a semi-classical  functional of class $s=1$.
The linear functional $\mathscr G_{\alpha,\beta}$ admits then the  integral representation
\begin{align*}
\big<\mathscr G_{\alpha,\beta},f\big>=\int_{-1}^{+1}\!f(x)|x|^{2\mu}(1-x^2)^{\alpha}dx,\ \forall f\in \mathscr P,
\end{align*}
and it is positive definite when $\alpha,\beta>-1$. Its  moments  are given by
 \begin{subequations}
 \begin{align}
&(\mathscr G_{\alpha,\beta})_{2n}=
\frac{\Gamma\big(\alpha+1\big)\Gamma\big(n+\beta+1\big)}{\Gamma\big(n+\alpha+\beta+2\big)}, \ n\geqslant0,\label{3.14a}\\
&(\mathscr G_{\alpha,\beta})_{2n+1}=0, \ n\geqslant0,\label{3.14b}
\end{align}
 \end{subequations}
 which are valid for  $\alpha+\beta\ne-n, n=1, 2, \ldots$, $\beta\ne-n, n=1, 2, \ldots$. In the above, the case $\beta\ne-1/2$ has to be excluded
 since this involves the ultraspherical or Gegenbauer polynomials ${G}_n^{\nu}$, with $\nu=\alpha+1/2$.\\
 Now, how can finding an orthogonalizing weight function for the polynomials ${G}_n^{(\alpha,\beta)}$ when $\alpha, \beta\ngtr-1$ on an interval
  of the real line. Here as well as in the previous case, the existence of such a weight is in fact guaranteed  by Boas's result \cite{Boas},
 and it may be constructed in any region in which $\alpha$, $\beta$ belong to the strips
 $-i-1<\alpha<-i$ and $-j-1<\beta<-j$, respectively,  where  $i$, $j$ are two positive integers.
 To do so, one can   proceed  in three independent  parts  to  regularize $w(x)=|x|^{2\mu}(1-x^2)^{\alpha}$ at $x=0$ and $x=\pm1$.\\
 Notice that the process requires a tedious computation in this case that we do not present here.\\
However, to have an idea about what we can expect when $\alpha, \beta\ngtr-1$, we refer the reader to \cite{MoKr}
where a detailed discussion concerning this technique was done, including  the Jacobi case.\\
\noindent We now return to apply  Theorem 2.2. First, a trivial verification shows that
$$
a_n=b_n=\frac{\xi_n}{n+\theta-1-\mu(-1)^n}, \,c_n=-\frac{2(\alpha+1)\xi_n}{n+\theta-1-\mu(-1)^n} \ \ \mbox{and}\ \,
 d_n=\xi_n\mu_n\frac{n+\delta_{n+1}+2\alpha+1}{n+\theta-1-\mu(-1)^n}.
 $$
 Next, when these are substituted in \eqref{eq2.17} we obtain after simplification that
\begin{align}
(1-x^2)T^2_\mu P_n-2(\alpha+1)xT_\mu P_n+\lambda_nP_n=0,\label{eq3.15}
\end{align}
where $\lambda_n=\mu_n\big(n+\delta_{n+1}+2\alpha+1\big),\ n\geqslant0,$  which  also writes
 \begin{align*}
\lambda_{2n}=4n\big(n+\alpha+\beta+1\big), \ n\geqslant0,\ \ \mbox{and}\ \
\lambda_{2n+1}=4\big(n+\alpha+1\big)\big(n+\beta +1\big), \ n\geqslant0.
\end{align*}
\noindent In this case, when $\mu=0$, we easily check that \eqref{eq3.15}, \eqref{eq3.13} are fully equivalent and deduce that
 they reduce to \eqref{eq2.31} with $a(x)=1-x^2$, $b(x)=-(2\nu+1)x$ and  $\lambda_n=n(n+2\nu)$, where $\nu$ stands for $\alpha+1/2$.\\
 Clearly, when $\mu$ assumes the value $0$, the three equations coincide and provide the well-known differential equation satisfied by the
  Gegenbauer polynomials ${G}_n^{\nu},\, n=0, 1, \ldots.$ \par
\medskip
\noindent {\bf Conclusion.}\\
 We studied the symmetric Dunkl-classical orthogonal polynomials in another way. This alternative approach for obtaining these polynomials
  has the advantage of being self-contained for determining the recurrence coefficients explicitly.
  A new characterization connecting the recurrence coefficients with the Dunkl-classical character is established in Theorem 2.1,
   while the second theorem (Theorem 2.2) provided another characterization of these polynomials by means of a second order differential-difference equation.
  Direct application of these two theorems  readily showed  that the only symmetric Dunkl-classical orthogonal polynomials are
  the generalized Hermite    and  generalized Gegenbauer  polynomials. These results are of course expected and consistent with those of \cite{BeGa}.

%\noindent {Acknowledgements}\\

\end{document}